\documentclass[oneside,leqno,a4paper,11pt]{article}
\usepackage{amssymb,amsmath,amstext,latexsym,theorem}
\usepackage{times}
\usepackage{epsfig}
\title{\huge \bf Locally homogeneous finitely nondegenerate CR-manifolds}
 
\author{\Large Gregor Fels }

\date{ }

\input \jobname.def

\begin{document}

\maketitle

\setlength{\abovedisplayskip}{10pt minus 2pt}
\setlength{\belowdisplayskip}{10pt minus 2pt}

\section{\llap{.\kern.75em}Introduction}\label{Introduction}

In several  areas of mathematics, homogeneous spaces are
fundamental objects as they
often serve as models for more general objects:  Various examples
from differential geometry (Riemannian symmetric spaces, principal bundles) and 
topology (geometric 3-manifolds), to algebraic and complex geometry 
(uniformization theorems, flag manifolds) 
etc.~underline the importance of spaces, furnished with a structure, compatible under a transitive group action. 
In this paper, we investigate homogeneous Cauchy-Riemann manifolds 
from the local point of view, more precisely, the germs of CR-manifolds which are
 locally homogeneous under some finite-dimensional Lie group.

\mn

The most common way of prescribing a  CR-manifold is to describe it locally 
in some $\C^n$
as the zero set  of certain defining functions. The characterization
of the geometric properties of such a manifold, like the signature of the Levi form(s),
finite or holomorphic nondegeneracy, minimality, etc. involves  a manipulation
of the defining equations, which, in concrete cases, can  be quite hard. 
A (locally) {\it homogeneous}\/ CR-manifold can also be 
described by a purely algebraic
datum, for instance by a CR-algebra in the sense of \cite{MN}. 
In fact, one can show that there is a natural equivalence between
the category of germs of locally homogeneous
CR-manifolds  on the complex geometric side 
and the category of 
CR-algebras on the algebraic side, see Section 4 for further details. 
In order to characterize the complex-geometric properties of $M,$  
the knowledge of the full Lie algebra
of local automorphisms of $M$ is not necessary; any Lie group, acting locally 
transitively on $M$ will do.
The advantage of this point of view is  that
in general the manipulation of CR-algebras is easier than the manipulation of the defining equations,
provided that there is a simple ``dictionary'' which ``translates'' the
algebraic properties of a given CR-algebra into the complex-geometric properties
of the underlying CR-manifold.

\mn  

In the first part of our paper we pursue  this goal and explain
how the Levi form of $M$ and its higher order analogues can be read off
 the corresponding CR-algebra. This enables us in Theorem \ref{LieOdege} 
to characterize the order of  nondegeneracy  of a locally homogeneous 
CR-manifold $M$,  as well as to decide whether or not $M$ is holomorphically degenerate.
In Theorem \ref{LieMini}, the minimality of $M$ is described in terms of the CR-algebra.
A basic ingredient in the proofs is the Main Lemma \ref{bas}, 
which relates certain canonical tensors and subbundles of $TM$ and $T^\C M$ 
with subspaces of infinitesimal CR-transformations and the corresponding 
Lie structure.
As a first application we generalize a result
 of Kaup and Zaitsev
stated in \cite{KZ2} (see the paragraph before \ref{kaza} for the precise statements)
for certain irreducible Hermitian symmetric spaces
to the more general case of arbitrary flag manifolds $Z$ with $b_2(Z)=1$ 
(Theorem \ref{kaza}). 
Our proof of this theorem  does not use Jordan-theoretical methods.

\mn

In the second part of this paper we provide an example of a homogeneous (hence, uniformly)
 3-nondegenerate hypersurface $\4M$ in the 7-dimensional Grassmannian of
isotropic 2-planes in $\C^7.$ In this example the first order 
Levi kernel is 3-dimensional
and contains the second order kernel  which is 1-dimensional.
While it is quite easy to produce 
real-analytic CR-manifolds which are, at some particular point,
 finitely nondegenerate of an arbitrary high order,
our hypersurface seems to be the first known example of a CR-manifold 
with a {\it uniform}\/ order of 
degeneracy bigger than 2. 
Note that in \cite{KZ2} orbits of real forms have been studied
in a certain  subclass  
of irreducible Hermitian symmetric spaces (of so-called tube type), and the authors prove
that all such orbits with a nontrivial CR-structure, (i.e., neither open nor totally real)
are 2-nondegenerate. Our example is an orbit in a more general flag manifold
and we use methods developed in the first part to determine its kind
 of nondegeneracy.
At this point one might expect to find orbits $M$ in complex flag manifolds 
$Z,$ with
uniformly finitely nondegenerate CR-structure of arbitrary high order, provided that the
ambient manifold $Z$ is general enough. Surprisingly,
at least for hypersurface orbits, this is not the case: In Theorem \ref{hypersurfOrb} we give a general upper bound for the order of degeneracy 
that is valid for all finitely nondegenerate hypersurface orbits 
in arbitrary flag manifolds. 
For instance, for all {\it classical}\/ cases, 
i.e., where the (connected component of the identity of the) complex group of biholomorphic transformations,
$\mathrm{Aut}(Z)^\circ,$ is a product of classical simple groups, this upper bound is  3. 
The methods used to determinate  the complex-geometric
properties of $\4M$ can be generalized to deal with 
arbitrary orbits of real forms in arbitrary flag manifolds.

\mn

We also like to mention that the methods developed in this paper
will be used in the forthcoming article \cite{FK2}, 
in which all 5-dimensional 2-nondegenerate germs of locally homogeneous CR-manifolds
are classified up to CR-equivalence.

\mn

Our paper is organized as follows. In Section 2 we discuss tensors induced
by Lie  brackets and, following Palais \cite{PAL},
 we recall basic facts on local actions. The main result here is the Main Lemma \ref{bas}.
Section 3 recalls basic geometric notions concerning CR-manifolds, focusing
on the condition of being finitely  nondegenerate.
In Section 4 we recall the definition of the category of CR-algebras (see also \cite{MN})
and show that there is an equivalence between this category and the category of germs of locally homogeneous CR-manifolds. The first part of our paper culminates in Section 5,
where we provide a  ``dictionary'', extracting from a given CR-algebra
the information necessary to characterize the complex-geometric properties of the 
underlying CR-germ. This characterization is used to prove a generalization of
the above mentioned result of Kaup and Zaitsev.
Finally, in Section 6 we give an example of a homogeneous 3-nondegenerate CR-mani\-fold 
(as already mentioned above) 
and  indicate a method how arbitrary orbits of real forms in flag manifolds can be handled.

\section{\llap{.\kern.75em}Tensors and homogeneous manifolds}

{\bf General notation.}
Let $X$ be a manifold. Given   a vector bundle  $\pi:\mathbb E\to X$  over $X,$
we write $\Gamma(X,\mathbb E)$ for the vector space of {\it smooth}\/ sections 
over $X$. If a further specification is necessary, we write 
$\Gamma^{\omega}\!(\cdot,\cdot)$ or
 $\Gamma^{\!\5O}\kern-1pt(\cdot,\cdot)$ etc.~for the real-analytic or
 holomorphic  sections, respectively.
By $\mathbb E_x$ we denote the fibre $\pi^{-1}(x)$ of $\mathbb E$ at $x\in X$.
As usual, $TX$ stands for the tangent bundle of $X$  and $T_xX$ for the tangent space at $x.$ 
Given a vector field $\xi\in \Gamma(X,TX)$ we write $\xi_x\in T_xX$ for 
its value at $x$.
If not otherwise stated all Lie groups and  Lie algebras 
(except for $\Gamma(X, TX)$) are assumed to be of finite dimension. In particular,
``homogeneous'' means (infinitesimally) homogeneous under a {\it finite}\/ dimensional Lie group (algebra). 
Lie groups are denoted by capital letters $G,H,..$ and the associated
 Lie algebras by the corresponding
fraktur letters $\7g,\7h,$ etc.
$G^\circ$ stands for the connected component of the identity of a Lie group $G.$
By definition, the Lie
bracket in $\7g$ is given by the Lie bracket of left-invariant vector fields on $G.$
By $\Ad$ we denote the adjoint representation of $G$ on $\7g$ and by $\ad$ its differential, i.e., $\ad_v(w)=[v,w].$
Given a real vector  space $V$, we denote by $V^\C:=V\otimes_\R\C=V\oplus iV$ 
 the formal complexification of $V.$ 
If the real vector space $V$ is furnished with an endomorphism $J:V\to V$ 
satisfying $J^2=-\Id,$ we write $V^{1,0}, V^{0,1}$ for the
$(\pm i)$\2--eigenspaces of $J^\C$ in  $V^\C.$

\mn

\bn
{\bf Tensors induced by Lie brackets.}
Let $\mathbb E\subset TX$ be a (smooth) 
subbundle.
It is well-known that the following 
$\R$\2--bilinear map 
\begin{equation}\label{ex1} \Gamma(X,\mathbb E)\times \Gamma(X,\mathbb E)\longrightarrow 
\QU{\Gamma(X,TX)}{\Gamma(X,\mathbb E)}\;,\quad (\xi,\eta)\mapsto [\xi,\eta]\kern.7em {\small\mbox{mod$\Gamma(X,\mathbb E)$}}
\end{equation}
is, in fact, $C^\infty(X)$\2--bilinear. Hence, it
induces a well-defined fibre-wise bilinear map (tensor) 
$\mathbb E_x\times \mathbb E_x\to \qu{T_xX}{\,\mathbb E_x}$, 
i.e., $[\xi,\eta]_x$mod$\mathbb E_x$
depends only on the values $\xi_x,\eta_x$ and not on the choice of the 
local sections $\xi,\eta$ in $\mathbb E.$

\mn

It turns out that for  (locally) {\it  homogeneous}\/
 manifolds $X$
the explicit computation of various tensors naturally attached to $X,$
 similar to that one given above,
can be reduced to a simple  algebraic expression.  
The main application we have in mind is the 
determination of the Levi form of a (locally) homogeneous CR-manifold $M$ and its ``higher-order'' analogues,
suitable for the characterization of the  
$k$\2--nondegeneracy of $M$ in the 
sense of  \cite{BHR}. In the next paragraphs we fix our notation
and briefly recall some basic
facts concerning homogeneous manifolds.  

\mn  

\bn
{\bf Locally homogeneous manifolds and bundles.} 
The topics of this subsection are well-known.
The reader familiar with the global concepts
 of a homogeneous space or a homogeneous bundle will have no difficulties to
give the local versions of these objects. 
In the following paragraphs we briefly recall the facts relevant for our purposes. A reference in the local situation is the fundamental paper of Palais 
\cite{PAL} (\cite{GO} is a more up-to-date reference).  

\mn

All groups occurring in this paper are assumed to be finite-dimensional
 Lie groups. Let $G$ be such a group. In the global
setting,  the fundamental objects are $G$\2--manifolds, i.e., manifolds
 provided with a (left) $G$\2--action $\cdot:G\times X\to X.$
A {\it homogeneous $G$\2-bundle} $\mathbb E\to X$ over such a manifold $X$
is a vector bundle together with a fibre-wise linear 
 action on $\mathbb E$ which is a lift
of the given  $G$\2-action on $X.$
If  $X$ is $G$\2--homogeneous, i.e., $G$ acts transitively on $X,$
we write $G_x$ for the isotropy subgroup at $x\in X$ and $\7g_x$ for the 
corresponding isotropy Lie subalgebra.
For a homogeneous bundle over a homogeneous manifold
 the isotropy representation
 $G_x\times \mathbb E_x \to \mathbb E_x$ determines
completely the global structure of the vector bundle $\mathbb E$
over $X=\qu G{G_x}\!:$ 
The total space of $\mathbb E$ is the twisted
$G$\2-product $G\tim{G_x} \mathbb E_x.$ 
Conversely,
a representation $H\to \GL(V)$ of a (closed) subgroup of $G$ on some vector space $V$
gives rise to the homogeneous vector bundle $\mathbb V:=G\tim{H}V$ over $G/H.$

\mn

All the above notions can be appropriately ``localized''.
A {\it local}\/ action of $G$ on a manifold $X$ is a map
$\cdot:\5U\to X$ such that $\5U\subset G\times M$ is an open neighbourhood of $\{e\}\times M$, the identity
$e\cd x=x$ holds for all $x\in X$ as well as $h\cdot (g\cdot x)=(hg)\cdot x$ when
 both sides are defined.
Without loss of generality we may assume that $G$ is simply connected, which
we do for all what follows.
A local action induces a Lie algebra homomorphism 
$\Xi:\7g\to \Gamma(X,TX),$ see  
\ref{ggtoTX}. A given  Lie algebra homomorphism $\Xi:\7g\to \Gamma(X,TX)$
is called  an {\it infinitesimal action}\/ of $\7g$ on $X$ 
and $X$ a $\7g$--{\it space}.
As shown in \cite{PAL},
an infinitesimal action $\Xi$ induces a local action of $G$ on $X$ 
(say, $G$ is the simply connected Lie group with Lie algebra $\7g$); consequently,
local and infinitesimal actions are equivalent objects.
It is known that  not globalizable local actions exist, see \cite{GO}, p.105 for further details.
All above notions can also be applied to germs of manifolds. 
To fix the notation, we write $(X,x)$ for a germ at the base point $x$ and $X$ for a representative
of the germ. Further, $(\Xi,X,x)$ stands for the germ  of a $\7g$\2--space 
(where the homomorphism $\Xi:\7g\to \Gamma(X,TX)$ describes the infinitesimal action).

\mn

By a {\it morphism}\/
between the $\7g$\2--space $X$ and the $\7g'$\2--space $X'$
we mean a pair $(\Psi,\psi),$ consisting of a map $\Psi:X\to X'$ in the given category (smooth, real-analytic, holomorphic)
and a Lie algebra homomorphism $\psi:\7g\to \7g'$ such that $\Psi_*(\Xi(v)_x)=\Xi'(\psi(v))_{\Psi(x)}$ for all $v\in \7g$. Every $\7g$\2--equivariant map (i.e., $\psi=\Id$) is an example of a morphism between two $\7g$\2--spaces.
\ A  morphism between the germs $(X,x), \ (Y,y)$ of two 
locally homogeneous spaces $X$ and $Y$
is then an equivalence class $[\Psi,\psi],$ induced by a  base point preserving
 equivariant morphism $(\Psi,\psi):X\to Y.$

\mn

We call an infinitesimal (or local) action of $\7g$ (resp. $G$) on $X$ 
{\it effective}\/ if the map $\Xi$ is injective. A global action 
$G\times X\to X$
is effective in this sense if and only if
the subgroup, formed by all elements  $g\in G$ which act as the identity on $X,$ is 
discrete.
Clearly, dividing $\7g$ (or $G$) by the
ineffectivity ideal $\7i=\ker \Xi$
(resp.~by the connected component of $\bigcap_{x\in X} G_x$), 
every non-effective action can be modified into an effective action with 
the same orbits 
(resp. Nagano-leafs).
A local, or equivalently, infinitesimal action on 
$X$ is called {\it transitive} if the evaluation map $\epsilon_x:\7g\to T_xX,$ $v\mapsto \Xi(v)_x,$
 is surjective for all $x\in X.$ 
We then say  that $X$ is {\it locally homogeneous}\/ or $\7g$\2--{\it homogeneous}. 
We call a germ $(X,x)$ {\it homogeneous}\/ if
there exists a locally homogeneous representative $X.$

\mn

It is known that for every pair $\7h\subset \7g$ of
 finite-dimensional Lie algebras,
there is a germ $(X,x)$ 
with a transitive  infinitesimal action $\Xi:\7g\to \Gamma(X)$ such that
 $\7h=\{v\in \7g:\Xi(v)_x=0\}=:\gx$. We call $\7g/\gx$
the {\it infinitesimal model}\/ for $(X,x).$ We say that the action or the infinitesimal model
 is effective if the action of $\7g$ 
 on some representative $X$ has this property.
In the case when $\7g$ is infinite dimensional, 
we do not know (even if $\dim \qu{\7g}{\gx}<\infty$) whether it is  always possible
to construct in a meaningful way a germ of a (finite dimensional) manifold with a local transitive action of some group ``associated'' with $\7g.$

\mn

Finally, a vector bundle $\pi:\mathbb E\to X$ over a $\7g$\2-homogeneous manifold 
is called {\it locally homogeneous} if the local action of $G$ lifts to a local action on $\mathbb E$ in  such a way that the corresponding local transformations are
fibre-wise linear. A germ of a locally homogeneous bundle 
(we use the notation $(\mathbb E,X,x)$ for it) 
is determined by the linear representation $\varrho:\gx\to \gl(\mathbb E_x)$
of the isotropy Lie algebra $\gx$ on the fibre $\mathbb E_x$.  
On the other hand, any representation $\varrho:\gx\to \gl(V)$ gives rise to a (germ of a) locally homogeneous
vector bundle $\mathbb V$ over the germ $(X,x)$ of a $\7g$\2--homogeneous manifold with $\mathbb V_x=V$.

\bn

Recall that each (local) $G$\2-action on $X$ 
induces the
so-called fundamental
vector fields on $X:$
The following map
\begin{equation}\label{ggtoTX}
\Xi:\7g\to \Gamma(X,TX),\quad v\mapsto \xi^v\;,
\mbox{\quad given by \  }\xi^v_y\,f:=\ddt \!\!\!f(\exp (-tv)\cd y) \;,
\end{equation}
where the $f\!$'s run through smooth functions defined in a neighborhood of $y,$
is a Lie algebra homomorphism. For each $v\in \7g$  the 
vector field  $\xi^v:=\Xi(v)$
 is called {\it fundamental}.
Unfortunately, the fundamental vector fields and (locally) homogeneous vector bundles
on a $\7g$\2--space seem to be unrelated.
For instance, the fundamental vector fields  
are {\it not}\/ invariant under the local group action.
Consequently, given a homogeneous $G$\2--subbundle $\mathbb E\subset TX$ 
and a fundamental vector field $\xi^v$ such that $\xi^v_x\in \mathbb E_x$ 
for some $x\in X,$  the values $\xi^v_y$ may not belong 
to $\mathbb E_y$ for $y$ close to $x.$  Since in general
the fundamental vector fields  do not generate a homogeneous subbundle 
$\mathbb E,$ they
  cannot be used  for the calculation 
of the Lie brackets in situations similar to \ref{ex1}.
Nevertheless the following lemma is valid, which is the main result of this section:

\begin{bale}\label{bas}
Let $X$ be a locally homogeneous $G$\2--manifold, $x\in X$ a base point
 and $\qu{\7g}{\7g_x}$ the corresponding infinitesimal model.
Let
$\mathbb E^1,\mathbb E^2,\mathbb D$ be any locally $G$\2-homogeneous 
subbundles of $TX.$ 
Let $\7e^1,\7e^2,\7d\subset \7g$ 
be the corresponding  $\ad(\7g_x)$\2--stable  
linear subspaces such that 
$\mathbb E^j_x=\7e^j/\7g_x$ and $\mathbb D_x=\7d/\7g_x.$ 
Assume that the bracket map 

\mn
\hfill 
$[\ ,\ ]:\Gamma(X,\mathbb E^1)\times \Gamma(X,\mathbb E^2)
\longrightarrow \qu{\Gamma(X,TX)}{\Gamma(X,\mathbb D)}$ \hfill\hbox{}

\mn is \ $C^\infty(X)$\2--bilinear,  
i.e., it defines a tensor 
$b:\mathbb E^1\oplus \mathbb E^2\to \qu{TX}{\kern1pt\mathbb D}$. 
For arbitrarily given tangent vectors 
$\nu^1\in \mathbb E^1_x,$  $\nu^2\in \mathbb E^2_x,$
choose  representatives $u^1\in \7e^1$ and $u^2\in \7e^2$.
Then, identifying $\qu{T_xX}{\mathbb D_x}$ with $\qu{\7g}{\7d}$, we have

\mn
\hfil $b_x(\nu^1,\nu^2)=[u^1,u^2]_\7g\kern4pt \mbox{\rm mod }\kern2pt \7d\;.$

\mn
Here, the bracket is taken in the Lie algebra $\7g,$
and the right-hand side does not depend on the choice of the representatives
$u^j.$
\end{bale}
\noindent
  
\Proof For simplicity, we carry out
the proof for globally homogeneous $X,$ i.e., $X=G/G_x,$ where  $G_x$ stands for the isotropy
subgroup at the base point $x.$ It relies on the construction of
particular local vector fields $\eta^1, \eta^2$ around $x$ and works equally well in the {\it locally}\/ homogeneous case.
By construction, the tensor $b$ is $G$\2--invariant. Hence, it suffices
to compute it at one point only.
Denote by  $\pi:G\to G/G_x$  the projection map and by 
$\pi_*:TG\to T(G/G_x)$ its 
differential. In particular, $\pi_*$ yields a surjection $\7g\to T_xX.$
Select once and for all a linear subspace $\7W\subset \7g,$ 
complementary to $\7g_x.$
 Let $v^1\in \mathbb E^1_x$ and  $v^2\in \mathbb E^2_x$ be arbitrarily given. 
Since $b$ is alternating, we may assume without loss of generality that
 $v^1,v^2$ are  linearly independent.
Select $w^1,w^2\in \7W$ such that $\pi_*(w^j)=v^j.$ Extend $w^1,w^2$
to a basis $w^1,\ldots,w^m$ of $\7W$ and let
$w^{m+1},\ldots,w^n$ be a basis of $\7g_x.$

\sn
By assumption, the bracket $[\eta^1,\eta^2]_x$mod$\mathbb D_x$
does not depend on the choice of the vector fields $\eta^j\in\Gamma(U,\mathbb E^j)$
which, for $j\in \{1,2\},$  extend $v^j$ in some neighborhood $U$ of $x.$
The key point here is the construction of  appropriate local extensions 
$\eta^1,\eta^2$ of $v^1 $ and $v^2.$ 
To accomplish this we first construct certain $\pi$\2--projectable vector 
fields $\zeta^j$ on an open set in $G$ and then 
define $\eta^j:=\pi_*(\zeta^j).$ 

\sn
{\sc Construction of the vector field $\zeta$ for a given $v\in T_xX.$}
Select a convex open neighbourhood $W\subset \7W\subset \7g$ of $0$ 
such that\hb
$\bullet$ the map $\exp :W\to  Y:=\exp(W)$ is a diffeomorphism onto the
locally closed submanifold $ Y\subset G$, and, 
\hb
$\bullet$ the restriction $\pi: Y \to X$ is a diffeomorphism onto
a neighborhood $V\subset G/G_x$ of $x,$ i.e., $ Y$ is the (image of a) 
local section in the principal bundle $\pi:G\to G/G_x.$
\hb
Write $(g,u)$ for elements in $TG=G\times \7g$ \wrt the trivialization
by left-invariant vector fields.
For an arbitrary given $v\in T_xX$ let $w\in \7W$ be the unique element with $\pi_*(w)=v.$
For such $w$  define $\zeta$ along $ Y$ simply by requiring 
$\zeta_g=(g,w)$ for all $g\in  Y$
and then extend $\zeta$ to a vector field 
 on $\pi^{-1}(V)= Y\cd G_x$ by
\begin{equation}\label{defze}
\zeta_{gh}:=(gh,\Ad_{h^{-1}} (w))\qquad g\in Y,\ h\in G_x \;.
\end{equation}
Note that $\zeta$ is invariant under the action of $G_x$ from the 
{\it right}; hence, it is $\pi$\2--projectable,
and we have $\zeta_{gh}=L^g_* R^h_* \zeta_e$, where $g\in  Y, h\in G_x.$
In particular, for the tangent vectors $v^1,v^2$ as above we write  $\zeta^1,\zeta^2$ for the
above constructed vector fields on $\pi^{-1}(V).$ 
Mutatis mutandis, this construction works
also in  the locally homogeneous situation.
 From the above follows that the vector fields
\begin{equation}\label{defxi}\eta^j:=\pi_*(\zeta^j),\quad j=1,2
\end{equation}
on $V\subset X$ are local sections in the $G$\2-bundles $\mathbb E^j$ with $\eta^j_x=v^j.$ (In general, the $\zeta$'s are neither
left- nor right-invariant.)
The $\pi$\2--projectable vector fields satisfy
\begin{equation}\label{proj}
 \pi_*[\zeta^1,\zeta^2]=[\pi_*\zeta^1,\pi_*\zeta^2]=[\eta^1,\eta^2]\;.
\end{equation}
We claim  that $[\zeta^1,\zeta^2]_e$
has a simple expression in terms of the Lie brackets in $\7g$ 
(by definition \wrt the  left-invariant vector fields).
Since $w^j, \ {1\le j\le n},$ form a basis of $\7g,$ the vector fields $\zeta^j$ 
can be written as linear combinations of left-invariant vector fields, i.e., 
 
\sn
\hfil$ \zeta^1=\sum_{j=1}^n a_jw^j_L\ ,\qqquad \zeta^2=\sum_{j=1}^n b_jw^j_L
$

\sn
with  $a_j,b_j\in C^\omega(\pi^{-1}V).$ 
By construction, all these functions are
constant on $Y$ and we have in particular $a_k\rest{Y}=0$ for $k\ne 1,$ and  $b_k\rest{Y}=0$ for $k\ne 2$.
The following identity is valid at an  arbitrary  point  
$y\in Y:$ 
$$
\begin{aligned}[]
[\zeta^1,\zeta^2]_y&=\textstyle\big[ \sum a_jw^j_L\ , \ \sum b_kw^k_L\big]=\cr
&=\textstyle\sum_{j,k} a_j(y)b_k(y)[w^j_L,w^k_L]+\sum_{j,k} a_j(y)(w^j_L b_k)w^k_L-\sum_{j,k} b_k(y)(w^k_La_j)w^j_L =\cr
&=\textstyle[w^1_L,w^2_L]+\sum_k (w^1_Lb_k)(y)\cd w^k_L-\sum_j (w^2_La_j)(y)\cd w^j_L 
\;.\cr
\end{aligned}
$$ 
Since $t\mapsto \exp tw^j\in Y$ are the local integral curves  at $e$ for $w^1_L$ and $w^2_L,$ 
it follows $w^1_Lb_k\ (e)=w^2_L a_k \ (e)=0$
for all $k,$
and the above formula, evaluated at $e,$ implies  $[\zeta^1,\zeta^2]_e=[w^1_L,w^2_L]_e.$ 
This identity together with \ref{proj} concludes our proof.
\qed
\begin{cor}\label{Jxi} {\bf (of the Proof of the Main Lemma)} 
\setlength{\parskip}{-5pt}
\begin{enumerate}
\setlength{\parskip}{-1pt}
\setlength{\partopsep}{-16pt}
\setlength{\itemsep}{+5pt}
\setlength{\itemindent}{0.0cm}
\setlength{\leftmargin}{12cm}
\item[$(i)$]Assume that $\mathbb E\subset TX$ is a
(locally) $G$\2--homogeneous vector subbundle over a (locally) homogeneous space $X$ and
$J:\mathbb E \to \mathbb E$ is a (locally) $G$\2--equivariant bundle endomorphism.
Let
$v\in \mathbb E_x$ be arbitrary and $v':=J_xv$. Select $w,w'\in \7W\subset \7g$ 
with $\pi_*(w)=v,$ $\pi_*(w')=v'$  and define $\zeta,\zeta'$ as in \ref{defze}.
 Then for the corresponding vector fields $\pi_*(\zeta)=\eta,\ \pi_*(\zeta')=\eta'\in \Gamma(V,\mathbb E)$ 
the relation $\eta'=J\eta$ holds at all points of $V\subset X.$
\item[\rm(\it ii\rm)]  
The statement of the Main Lemma remains true if $TX$ is replaced by its formal complexification 
$T^\C X=TX\otimes_\R \C=G\tim{H}\qu{\gc}{\gc_x}$ and $\mathbb E^1,\mathbb E^2,\mathbb D$ are $G$\2-homogeneous subbundles, corresponding  to the linear subspaces $\7e^1,\7e^2,\7d$ of $\gc$. Further, the Main Lemma remains true if the
tensor $b$ is defined by a linear combination of brackets 
(even if every  single bracket, which occurs in such an expression, does not
 yield 
a well-defined tensor).

\end{enumerate}
\end{cor}
 
\mn
In the next section we apply the formula stated in the Main Lemma 
to locally homogeneous CR-manifolds
for the computation of their Levi forms and certain higher order analogues. 
This will enable us 
to give a simple characterization of the (non)degen\-eracy type 
for  locally homogeneous  CR-manifolds.

\section{\llap{.\kern.75em}CR-manifolds and nondegeneracy conditions}
\setcounter{equation}{0}
\label{CRgeom}

In this section we briefly recall some
basic facts concerning CR-manifolds and certain geometric properties of them.
In particular, we closely examine the condition of being finitely 
nondegenerate, which 
plays a major role in the next sections.
As a general reference for CR-manifolds, 
see \cite{BER}
and \cite{Bog}.

\begin{defi}
\label{CR-Mf}
An abstract CR-manifold  is a smooth manifold $M$ together with a
  subbundle $\H\subset TM$ (we call it the {\it complex subbundle}) 
and a vector bundle
 endomorphism $J:\H\to \H$ with $J^2=-\Id$ 
(the so-called {\em partial almost complex structure}) such that 
for all $\xi,\eta\in \Gamma(X,\H)$ 
it follows\footnote[1]{In \cite{KZ1}, Sec.~2 an even more general definition of a CR-manifold has been given. }  $[\xi,\eta]-[J\xi,J\eta] \in \Gamma(M,\H).$ \  
If, in addition, the Nijenhuis tensor 

\mn
\hfill $N(J)(\xi,\eta)=[J\xi,J\eta]-[\xi,\eta]-J([\xi,J\eta]-[J\xi, \eta]), $ $\xi,\eta\in \Gamma(M,\mathbb H)$, \hfill{}

\mn
  of $J$ vanishes,  
we call $(M,\H,J)$ {\it formally integrable}.
\end{defi}
In this paper we almost exclusively investigate manifolds which are
locally homogeneous under some Lie group.
Every smooth manifold
furnished with a smooth locally transitive action of a 
finite dimensional Lie group automatically carries a real-analytic structure,
 compatible with the group action. 
We assume from now on (if the contrary is not explicitly stated) that {\bf all manifolds, actions and subbundles are real-analytic and the CR-manifolds are formally integrable.}
However, the sections in such subbundles may be only  smooth.

\mn

Two "extreme" classes  of  CR-manifolds are the following: 
Complex manifolds $Z$ are 
 precisely those formally integrable CR-manifolds with maximal possible 
complex subbundle: $\mathbb H=TZ.$ Here, $J:TZ\to TZ$ 
is the complex structure,
induced by the multiplication with $i=\sqrt{-1}$ in local coordinate charts.
On the other hand,
every real manifold, furnished  with the trivial CR-structure $\H=0$ is CR and called
{\it totally real}\/ as a CR-manifold.

From the local point of view 
complex manifolds as well as real manifolds with  
$\H=0$ are not very interesting. Hence,
apart from few exceptions, the  CR-manifolds considered in
this paper do not belong to any of the above two classes.
A wide class of CR-manifolds consists of
 real submanifolds $M$ of complex manifolds 
$(Z,J)$ such that $\H_x:=T_xM\cap JT_xM$ and $\dim \H_x$ is a constant function of $x\in X.$ 
Such a CR-manifold is formally integrable (since $(Z,J)$ has this property).
On the other hand, due to the well-known embedding theorem of 
Andreotti-Fredricks (\cite{AF}),
every formally integrable real-analytic CR-manifold
admits a generic CR-embedding into a complex manifold $Z.$ 
Hence, without loss of generality we assume in the following that all 
CR-manifolds
under considaration  are (locally) closed submanifolds $M\into Z$
and fulfill the above conditions together with  $TZ|_M=TM+ JTM$ (genericity).

\bn
{\bf Infinitesimal CR-transformations.}
Let $M=(M,\H,J)$ be a real-analytic CR-manifold.  
There is a particular Lie subalgebra of $\Gamma(M,TM)$, 
related to the CR-structure:
 Call a vector field $\xi\in \Gamma^\omega(M,TM)$ 
an {\it infinitesimal CR-transformation}\/
 if the corresponding  local 1-parameter subgroup $\Psi^\xi_t$ of $\xi$
acts by local CR-transformations of $M$. Write $(M,o)$ for the germ at $o\in M$
of $M.$ Define
$\hol(M)\subset \Gamma^\omega(M,TM)$ (resp.~$\hol(M,o),$ if dealing with germs) 
as the subspace consisting of (germs of) infinitesimal CR-transformations
of $M$ (or $(M,o)$, respectively; the elements in   $\hol(M,o)$ \nn  
vanish at $o$).
The spaces $\hol(M)$ and $\hol(M,o)$ are
Lie algebras, with Lie structure induced
 by the usual Lie brackets of vector fields. 
In the above definition we do {\it not
 require}\/ that
the infinitesimal CR-transformation on an embedded CR-manifold, $M\into Z,$
are restriction of holomorphic vector fields on $Z.$ However,
due to  Proposition 12.4.22 in \cite{BER}, this follows automatically.
Finally, by a 
{\it holomorphic vector field}\/ on a complex manifold $Z$ we mean
a holomorphic section in the real tangent bundle $TZ.$
Given a manifold $M$ with some structure $\7C,$ we write $\mathrm{Aut}_\7C(M),$ or simply $\mathrm{Aut}(M)$ for the group of 
all automorphisms of $M$ preserving this  structure and $\mathfrak{aut}(M)$ for the correspodning Lie algebra.

\bn
{\bf The notion of \9k-nondegeneracy.}
A basic invariant of a CR-manifold is its vector-valued Levi form 
$\4L^M$, or equivalently \wrt the encoded information,
the canonical alternating 2-form
$\omega^M:\H\oplus \H\to \qu{TM}{\H}.$ 
This 2-form  is simply  the tensor induced by Lie brackets (as in \ref{ex1}).
The (classical) Levi form\footnote[2]{We took the definition from \cite{KZ1}. 
It differs from the Levi form considered by some other authors, see for instance \cite{BER}, by the factor $i/2$} $\4L^M,$ 
which is a $J$\2--invariant sesquilinear tensor 

\sn
\hfil
$\4L^M:\H\oplus \H \to \qu{T^\C M}{\H^\C},$ \hfill

\sn
 and $\omega^M$ are related:
$\4L^M(u,v)=\omega^M(u,v)+i\omega^M(Ju,v).$
A complexified version of the Levi form is the tensor $\4L^1:\Ha\oplus \Hh\to \qu{T^\C M}{\kern1pt \H^\C}$ induced by Lie brackets of local sections in $\Ha$ and $\Hh.$
Set
\[
\Fa_{(0)}:=\Ha_{\hbox{ }},\qquad 
\Fa_{(1)}:=\{\xi\in \Ha:\4L^1(\xi\ ,\ \Hh)=0\}.
\]
A CR-manifold is called {\it Levi-nondegenerate}\/ or 1--{\it nondegenerate at} $x\in M$ 
if the fibre of  $\Fa_{(1)}$ at $x$ is  zero.

The notion of $k$\2-nondegeneracy of $M$ at a point $x$ 
has  been originally  defined in \cite{BHR} (see also Sec.~11.1 in \cite{BER}) for arbitrary CR-manifolds.
In general, the order $k$ of nondegeneracy at $x\in M$ varies from point to point and can be arbitrarily high.
For the class of CR-manifolds of ``uniform degeneracy'' 
(i.e., the dimensions  of  all  fibre-wise defined subspaces 
$(\Fa_{(k)})^{\hbox{\ }}_x\subset T^\C_xM$, as constructed below,
 do not depend on  $x\in M$ and  form 
well-defined sub{\it bundles} of $T^\C M$)
which includes all locally homogeneous CR-manifolds,
$k$\2-nondegeneracy can be expressed as 
the nondegeneracy of certain tensors $\4L^{k+1}.$ The latter tensors
can be considered as a generalization of the Levi form 
$\4L^1.$ This has already been explained in the Appendix of \cite{KZ2}. 
For convenience, we recall this construction in a form suitable for our purposes.

Define recursively the subbundles
\begin{equation}
\label{Fk}
 \Fa_{(k)}:=\{\xi\in \Fa_{(k-1)}:\4L^{k}(\xi\ , \ \Hh)=0\}\;,
\end{equation}
and the following maps, induced by Lie brackets:
\begin{equation}
\label{Lk}
 \4L^{k+1}:\Fa_{(k)}\times \Hh_{\raise5pt\hbox{ }} \longrightarrow \QU{\Fa_{(k-1)}\oplus \Hh_{\raise5pt\hbox{ }}\!\!}{\Fa_{(k)}\!\oplus \Hh_{\raise5pt\hbox{ }}} \quad \subset \quad \QU{\H^\C\!}{\Fa_{(k)}\!\oplus \Hh_{\raise5pt\hbox{ }}\;.} 
\end{equation}
The fact that all $\4L^{k+1}$'s are  well-defined tensors follows from the formula 
$-d\theta(\ph ,\eta)=\theta([\ph,\eta])$, where $\ph$ and $\eta$ are  local sections in 
$\Fa_{(k)}$ and $\Hh,$ 
respectively, and the $\theta$'s run through all
 1-forms $\theta:T^\C M \to \C$ which vanish on 
$\Fa_{(k)}\oplus \Hh_{\raise5pt\hbox{ }}.$
By construction, for each CR-manifold of uniform degeneracy there is 
the following filtration of $\Ha$ by complex subbundles:
 $\Ha=\Fa_{(0)}\supset \Fa_{(1)} \supset \Fa_{(2)} \supset \cdots .$
The property of being 
$k$\2--nonde\-generate is characterized  in the following

\begin{prop} 
\label{knondeg}
Let $\Fa_{(j)},\ j=0,1,2,...$ be the subbundles as defined in \ref{Fk}.
A 
 CR-manifold $M$ of uniform degeneracy 
is  $k$\2--non\-degenerate if and only if 
$ \Fa_{(k-1)}\ne\Fa_{(k)}=0.$  
\end{prop}

For locally homogeneous CR-manifolds the subbundles
and tensors, as defined in \ref{Fk} and \ref{Lk}, respectively,
 can be characterized in Lie algebraic terms.
In particular, the geometric notion of $k$\2--nondegeneracy can be completely 
described in terms of a filtration of certain subalgebras, 
as  will be shown in Section 5.

\section{\llap{.\kern.7em}Homogeneous CR-germs and CR-algebras}
\label{CRgerms}
\setcounter{equation}{0}

In this section we show that each germ $(M,o)$
of a locally homogeneous real-analytic CR-manifold 
(homogeneous CR-germ, for short)
can be described by an algebraic datum, for instance by a CR-algebra.
Vice versa, every CR-algebra gives rise to a homogeneous CR-germ and all these
 assignments are functorial. 
We start by recalling the definition of the 
category of CR-algebras, essentially following \cite{MN}.

\bn
{\bf The category of  CR-algebras.}
To fix notation, let $\7g$ stand for a real Lie algebra,
 let $\gc:=\7g\otimes_\R \C=\7g\oplus i\7g$
be its complexification and $\psi^\C$  the complexification
of a real homomorphism $\psi:\7g\to \7g'.$ As before, we write $\7l$ for the 
complexification $\gc$ and $\sigma$ for the unique complex conjugation 
$\7l\to \7l,$ fixing the real form  $\7g\subset \7l.$

\mn

A pair, consisting of a finite-dimensional real Lie algebra $\7g$
and a complex subalgebra $\7q$ of $\7l:=\gc$ is called 
a {\it CR-algebra}. In contrast
to \cite{MN}, here we require the finite dimensionality of $\7g.$ 
A {\it morphism} $(\7g,\7q)\to (\7g',\7q')$ is a Lie algebra homomorphism
$\psi:\7g\to \7g'$ with $\psi^\C(\7q)\subset \7q'.$
We refer to the category in which the objects are CR-algebras and the morphisms are as just described
as to the category of {\it CR-algebras}, or, for short,  \AHCR\kern-2pt.

\bn

On the geometric side there is the 
{\bf category of  homogeneous CR-germs.}
The objects in this category
are homogeneous CR-germs $(\Xi,M,o)$ and the morphisms $[\Psi,\psi]$
are 
as defined 
in the subsection ``{\sl locally homogeneous manifolds and bundles}'' of section 2.
Note that $\Psi$ automatically is a CR-map.
We refer to this category as to the category of 
{\it homogeneous CR-germs} (and write \GHCR, for short).

\sn

Discarding for a moment local actions,
 there is also the category  \CRo, consisting of germs of real-analytic
CR-manifolds as objects and real-analytic (germs of) base point preserving 
CR-maps $(M,o)\to (M',o')$ as
morphisms.
We have then the obvious forgetful functor  \GHCR $\rightsquigarrow $ 
\CRo.
Note, however, that the notion of an isomorphism 
is different in these two categories: Two homogeneous CR-germs $(\Xi,M,o)$ and 
$(\Xi',M',o')$ may be non-isomorphic in \GHCR, 
though the underlying CR-germs are
CR-equi\-valent, i.e., isomorphic in \CRo.
To distinguish these two notions of an isomorphism, 
we refer to  $(\Xi,M,o)$ and $(\Xi',M',o')$ as {\it isomorphic}\/
if there is an isomorphism between them in \GHCR and as {\it CR-equivalent}\/
if $(M,o)$ and $(M',o')$ are isomorphic in \CRo.
This fine point plays a role in \cite{FK2}, where 5-dimensional 2-nondegenerate
homogeneous CR-germs are classified up to CR-equivalence, 
and this classification is reduced to the classification of 
$\7g$\2--homogeneous CR-germs with $\dim \7g$ as small as possible.

\bn
{\bf Functors.}
There is a functor $\4G$ from the category of CR-algebras to the category of
homogeneous CR-germs (this has also been remarked in \cite{MN}).
Given a CR-algebra $(\7g,\7q)$ set $\7l:=\gc.$
Let $(Z,o)$ be the germ of a complex homogeneous manifold with the 
infinitesimal model $\qu{\7l}{\7q}$ and $Z$ a locally homogeneous representative.
 The CR-germ $(\Xi, M,o)$ is then determined
as the germ at $o$ of the real-analytic Nagano leaf $M$ through $o$ in 
$Z$ \wrt $\7g\otimes_\R^{\hbox{\ }}C^\omega(Z)$ (see \cite{Na}).

Let $\psi:(\7g,\7q)\to (\7g',\7q')$  be a morphism between two CR-algebras.
Let $Z$ and $Z'$ be  representatives of the germs of complex manifolds, 
determined by the 
infinitesimal models $\7l/\7q$  and $\7l'/\7q',$ respectively.
Then $\psi$ induces 
(possibly after  shrinking  $Z$) an $(\7l,\7l')$\2--equivariant, 
holomorphic and
 base point preserving map $\Psi:Z\to Z',$
which maps $M\subset Z$ to $M'\subset Z'.$ Hence, the restriction of $\Psi$
to $M$ is a real-analytic CR-map and yields  a morphism
  between 
the homogeneous CR-germs $(\Xi, M,o)$ and $(\Xi',M',o').$

\mn 

There exists also a functor $\4A$  in the  opposite direction. 
Let a $\7g$\2--homogeneous CR-germ $(\Xi,M,o)$ be given. 
Due to \cite{AF}, there exists a
 complex manifold $Z$ such that a representative $M$ 
is generically CR-embedded in $Z.$ 
The only point here is that this embedding is automatically 
locally equivariant \wrt $\7g:$ This
is a consequence of the extension results in \cite{BER} (Corollaries 12.4.17 and 1.7.13)
and our assumption that $\7g$ is finite dimensional. Hence, 
possibly after shrinking  $Z,$
we may assume that for each $v\in \7g$ the vector field $\Xi(v)$ is the restriction of a 
holomorphic  vector field on $Z$. 
Therefore, we can consider $\Xi$ as a Lie algebra 
homomorphism $\7g\to \Gamma^{\!\5O}(Z,TZ).$
Since the  Lie algebra $\Gamma^{\!\5O}(Z,TZ)$ is complex, $\Xi$ extends to a complex homomorphism
$\Xi^\C:\7l \to \Gamma^{\!\5O}(Z,TZ).$ Define the complex isotropy subalgebra
$\7q:=\{w\in \7l:\Xi(w)_o=0\}.$
The pair $(\7g,\7q)=:\4A(\Xi, M,o)$ is a CR-algebra and we call it the CR-algebra {\it associated with} $(\Xi, M,o).$  
Define $\go:=\7g\cap \7q.$ Observe that $\7g/\go$ is the infinitesimal model for 
$(\Xi,M,o)$ and
$\7l/\7q$ the infinitesimal model for $(\Xi^\C, Z,o).$
A word of caution: Even if $\Xi:\7g\to \Gamma(Z,TZ)$ is injective, i.e., the original
$\7g$\2--action is effective, the complexification $\Xi^\C$ may not be injective, i.e., the sum
$\Xi(\7g)+J\Xi(\7g)$ may not be direct.

It follows that an equivariant morphism 
$(\Psi,\psi):(\Xi,M,o) \to (\Xi',M',o')$
induces a morphism of the associated CR-algebras: The only point which has to be checked is
that the complexification of $\psi:\7g\to \7g'$ maps $\7q$ to $\7q':$
Again by the extension results from \cite{BER}, a representative $\Psi:M\to M'$ extends
to a holomorphic map $\widehat\Psi:Z\to Z'.$ By the identity principle, $\widehat\Psi$ is
equivariant \wrt $\7l$ and $\7l'.$ Since $\widehat\Psi$ preserves the base points, the inclusion $\psi^\C(\7q)\subset \7q'$ follows from 
$\widehat\Psi_*(\Xi(w)_o)=\Xi'(\psi^\C(w))_{\hat\Psi(o)}=\Xi'(\psi^\C(w))_{o'} .
$
Summarizing, we have

\begin{prop}
\label{IsoOfCat}
The above defined covariant functors 

\mn
\centerline{\AHCR $\stackrel{\4G}{\rightsquigarrow}$ \GHCR
 \quad and \quad 
\GHCR $\stackrel{\4A}{\rightsquigarrow}$  \AHCR }

\mn
are  mutually quasi-inverse and
yield an  equivalence of the two categories.
\end{prop}

\mn
{\bf Remark.} There  exist (locally) homogeneous manifolds with non-integrable CR-structures. A
germ of such a more general CR-manifold can also be 
described by purely algebraic data, 
for instance by a  quadruple $(\go,\7g,\7H,J),$ consisting of the 
Lie algebras $\go\subset \7g,$   an $\ad(\go)$\2--stable subspace $\7H$ of $\7g$
and an endomorphism $J:\qu{\7H}{\go} \to \qu{\7H}{\go}$ such that
$J$ is $\ad(\go)$\2--equivariant, $J^2=-\Id$ and
$[\widehat Jv,\widehat Jw]-[v,w]\in \7H$ holds for all $v,w\in \7H$ and some linear lift $\widehat J:\7H\to \7H$ of  $J$ with $\widehat J(\go)\subset \go.$
 However, such quadruples $(\go,\7g,\7H,J)$ seem to be
less convenient to deal with than  CR-algebras.

\bn

\section{\llap{.\kern.75em}Geometric properties of a germ, given by a CR-algebra}
\setcounter{equation}{0}

As seen in the previous section, 
the germ at $o$ of a locally homogeneous CR-manifold $M$ is completely determined
by the corresponding CR-algebra.
Consequently,
 all objects naturally attached 
to $M$ and their geometric properties are (at least a priori) 
completely determined by $(\7g,\7q).$ 
In this section we show in an explicit way how the geometric 
information encoded in a CR-algebra
can be extracted. In particular, we give a description of the subbundles 
$\H,\Ha,\Hh,\Fa_{(k)}$ of $T^\C M$ in terms of quotients of Lie algebras.
The main results of this section are  a description of the $k$\2--nondegeneracy
and the holomorphic nondegeneracy of a CR-germ $(M,o)$ as a purely algebraic property
of its CR-algebra (Theorem \ref{LieOdege}; see also the following remarks), 
and  Theorem \ref{LieMini} in which the minimality of $M$ 
is characterized
in a similar fashion.  
As an application we give a simple proof of the following result: 
each  non-extreme $G$\2--orbit in $Z=L/Q,$ 
where $Z$ is an arbitrary  flag manifold with $b_2(Z)=1,$ $L$ a complex subgroup of $\mathrm{Aut}(Z)$  and $G\subset L$ an arbitrary real form,
is minimal and holomorphically nondegenerate. This generalizes a theorem of 
Kaup and Zaitsev, see \cite{KZ2}.

\bn

Let $(\7g,\7q)$ be a given CR-algebra.
Recall that  $\7l=\gc,$ $\go:=\7g\cap \7q$
and $\sigma:\7l\to \7l$ is the involutive automorphism with
$\7l^\sigma= \{v\in \7l:\sigma(v)=v\} =\7g.$ 
Let $(M,o)$ be the corresponding homogeneous CR-germ which is CR-embedded in $(Z,o),$ 
as explained in section \ref{CRgerms}.
Since the  vector bundles $TM, \H, \Hh, \Ha, T^{1,0}Z,$ etc.~are locally homogeneous
\wrt the given transitive local actions on $M$ and $Z$, they are 
determined by a single fibre, say at
$o\in M.$  As these various
fibres are subspaces of the corresponding (complexifications of) tangent spaces
$T_oM=\qu{\7g}{\go},$ 
$T_oZ=\qu{\7l}{\7q},$  $T^\C_oM=\qu{\gc}{\7g^\C_o}$ etc., we need to specify
the appropriate subspaces of the preceeding quotients of Lie algebras.
We proceed with  preparatory observations.

\sn
$\bullet$ The real isotropy Lie algebra $\go$ is a real form
of $\7q\cap \sigma\7q$ (this was already observed in \cite{W}). Hence,
the complexified tangent space $T^\C_o M$ is the quotient $\qu{\7l}{\qq}.$

\noindent
$\bullet$
Define the subspace $\7H:=(\7q+\sigma\7q)^\sigma=(\7q+\sigma\7q)\cap \7g$ of $\7g.$ Note that $[\go,\7H]\subset \7H$ and  observe that the map 
$\7q\to \7H,$ $w\mapsto w+\sigma w$ is surjective. The quotient
$\qu{\7H}{\go}$ coincides with the intersection $\qu{\7g}{\go} \cap i(\qu{\7g}{\go}):$
This follows from the equation $\7H=\{x\in \7g: x+\7q=iy+\7q$ for some $y\in\7g\}.$

\noindent
$\bullet$
The invariant complex structure $J:\H\to \H$ induced by the embedding $M\into Z,$ 
i.e., the endomorphism $J_o:\qu{\7H}{\go} \to \qu{\7H}{\go},$ can be described
 as follows:
Recall that given
any $u\in \7H$ there exists a $w\in \7q$ with $u=w+\sigma w.$
Further, since $\7l=\7g\oplus i\7g,$ each element in $\7l$ has the unique
decomposition into its real and  imaginary parts. Then:  
\begin{equation}
\label{JfromEmbe} 
J_o:\qu{\7H}{\go} \to \qu{\7H}{\go},\qquad (u+\go)\mapsto 2y+\go \;,
\end{equation}
where $u=w+\sigma w$ mod $\go,$  $w\in \7q$ and $x+iy$ is the decomposition of $w$ into its
 real and  imaginary parts.

We summarize the above results, i.e., the identifications of the various fibres at $o$ with the corresponding quotients of Lie algebras in the following diagram:
\begin{equation}
\label{Diagra}
\diagram{\qu{\7H}{\go} & = &\H_o \cr
\cap & & \cap \cr
\qu{\7g}{\go} & \ = \ & \tom & \rinto & \toz & \ = \ & \qu{\7l}{\7q} \cr
\cap & & \cap & & \cap & &\cap \cr
\qqu {\7l}{\qq} & = & T^\C_{o}M & \rinto & T^\C_{o}Z & = & \qqu {\7l^\C}{\qc} \cr
\cup & & \cup & & \Vert & &\Vert \cr
\qqu {\sigma\7q}{\qq} & = & \Hho & \rinto & T^{1,0}_{o}Z & = & \qqu{\7l^{1,0}}{\7q^{1,0}}  &  \ = \qqu {\7l}{\7q}\cr
\oplus & &\oplus & &\oplus & &\oplus \cr
\qqu {\7q}{\qq} & = & \Hao & \rinto & T^{0,1}_{o}Z & = & \qqu{\7l^{0,1}}{\7q^{0,1}}&  \ \buildrel\sigma\over= \qqu {\7l       }{\sigma\7q}
}
\end{equation}

\bigbreak
\ 

\bigbreak
\noindent
{\bf Finite nondegeneracy in terms of  CR-algebras.}
In the next paragraphs  we  repeatedly  apply  the Main Lemma \ref{bas} to
the various tensors associated with a locally homogeneous CR-manifold $M$ as
 described in  Section \ref{CRgeom}. We 
obtain in that way expressions for all $\4L^k$'s 
in terms of Lie brackets in the 
Lie algebra $\7l.$ Here, $\7l=\gc$ comes from the CR-algebra, associated to a 
given $\7g$\2--homogeneous CR-germ $(M,o).$  
Keeping in mind the identifications \ref{Diagra}, 
the Main Lemma \ref{bas} immediately implies
\begin{equation}
 \omega^M:\qu{\7H}{\go} \times \qu{\7H}{\go}\to \qu{\7g}{\7H}, \qquad (u,v)\mapsto [u,v]_\7g\hbox{ mod }\7H \;. 
\end{equation}
As already mentioned, the complexification of $\omega^M,$ restricted to $\Ha\times \Hh,$
i.e.,  the invariant tensor $\4L^1,$
is equal to the  Levi form up to some  factor. 
Also in this case the Main Lemma together with the identifications \ref{Diagra} implies the following formula
for  $\4L^1$  at $o:$ 
\begin{equation}
 \4L^1: \qu{\7q}{\qq} \times  \qu{\sigma\7q}{\qq} \to \qu{\7l}{(\7q{+}\sigma\7q)}\ ,
(u,v)\mapsto [u,v]_\7l \hbox{ mod } (\7q{+}\sigma\7q)\;.
\end{equation}
For short, write $\7q^{(0)}:=\7q,$ $\7q^{(\infty)}:=\7q\cap\sigma\7q.$ The (left-) kernel of $\4L^1$ is $\qu{\7q^{(1)}}{\qq},$ 
where 
\begin{equation}
\label{f1}
\7q^{(1)}:=\{ w\in\7q: [w,\sigma\7q]\subset \7q+\sigma\7q\}
\end{equation} 
coincides with the normalizer $N_{\7q}(\7q{+}\sigma\7q)$ and consequently
 is a complex subalgebra.
Similarly, the recursively defined (\ref{Lk}) tensors $\4L^k$
(which are invariant under the local action)  
are given by the formulae:
\begin{eqnarray}
\label{LieLk}
\4L^{k+1}:\qu{\7q^{(k)}}{\7q^{(\infty)}}&\times & \qu{\sigma\7q^{(0)}}{\7q^{(\infty)}}\longrightarrow
\QU{\7q^{(k-1)}+\sigma\7q^{(0)}\ }{\7q^{(k)}+ \sigma\7q^{(0)}} \\
 (u&,&v)\kern3.2em \longmapsto \quad [u,v]\hbox{ \ \small mod $\7q^{(k)}\!+\sigma\7q^{(0)}$}.\nonumber 
\end{eqnarray}  
Here and above, the right-hand sides does not depend on the choice of the
representatives $u$ and $v.$
The (left) kernels of $\4L^{k+1}$ are the {\it homogeneous}\/
 subbundles $\Fa_{(k+1)}$ (see \ref{Fk}); hence, 
they are   determined by the corresponding fibres at $o.$ 
A glance at (\ref{LieLk}) suggests the following definition:
\begin{equation}\label{ffk}
\7q^{(k+1)}:=\{w\in \7q^{(k)}:[w,\sigma\7q^{(0)}]\subset \7q^{(k)}+\sigma\7q^{(0)}\}\;.
\end{equation}

\begin{obser}
\label{Fkato}
The fibre of $\Fa_{(k)}$ at $o$ is isomorphic to the quotient $\qu{\7q^{(k)}}{\7q^{(\infty)}}.$
\end{obser}
Next, we prove  the auxiliary
\begin{lemma}
\label{ffilt} Let $(\7g,\7q)$ be a CR-algebra, $\7H=(\7q+\sigma\7q)^\sigma$ 
 and let $\7q^{(k)}$ be the
subspaces
of $\7q,$ defined in \ref{ffk}. Then
\setlength{\parskip}{-6pt}
\setlength{\partopsep}{-6pt}
\begin{enumerate}
\setlength{\itemsep}{-1pt}
\setlength{\leftmargin}{-2cm}
\item[(i)]
The real subspace $\7F:=N_\7g(\7H)\cap \7H$ is a subalgebra
and 
$(\7q^{(1)}+\sigma\7q^{(1)})^\sigma=\7F.$ 
\item[(ii)]
All subspaces occurring in the filtration 
$\7q=\7q^{(0)}\supset \7q^{(1)}\supset \7q^{(2)}\supset \ \cdots\  \supset\7q^{(\infty)} $
are complex subalgebras of $\7q.$

\sn

\end{enumerate}
\end{lemma} 
\Proof
ad (i): To show that $\7F$ is a subalgebra, it suffices to show that for $u,v\in \7F,$
$[u,v]$ belongs to $\7H.$ This  follows from $[\7F,\7F]\subset [\7F,\7H]\subset \7H.$
For the proof of the second identity note that the inclusion
$\7q^{(1)}+\sigma\7q^{(1)}\subset N_\7l(\7q{+}\sigma\7q)\cap (\7q{+}\sigma\7q)=(\7F^{(1)})^\C$
is obvious. Let now $u+\sigma w\in N_\7l(\7q{+}\sigma\7q)\cap (\7q{+}\sigma\7q)$
 be an arbitrary element with $u,w\in \7q.$ 
If $[u,\7q{+}\sigma\7q]\not\subset \7q{+}\sigma\7q,$ i.e., if there were $a\in \7q$ with 
$[u,\sigma a]\ne 0$ mod $\7q{+}\sigma\7q,$ 
then also $[u+\sigma w,\sigma a]\ne 0$ mod $\7q{+}\sigma\7q,$ contrary to the construction
 of $u+\sigma w.$ It follows $u,w\in \7q^{(1)}.$
\hb
ad (ii): Clearly, $\7q=\7q^{(0)}$ and $\7q^{(1)}$ are  subalgebras.
Assume  that  we have already proven that $\7q^{(j)}$ are subalgebras for
all $j<k.$ To conclude that $\7q^{(k)}\subset \7q^{(k-1)}$ is also a subalgebra, 
note that for $u,v\in \7q^{(k)}$ we have 
\setlength{\abovedisplayskip}{6pt minus 2pt}
\setlength{\belowdisplayskip}{6pt minus 2pt}
\begin{eqnarray*}
[[u,v],\sigma\7q]&\subset &[u,[v,\sigma\7q]]+[v,[u,\sigma\7q]]\subset [u,\7q^{(k-1)}]+[v,\7q^{(k-1)}]+\7q^{(k-1)}+\sigma\7q\subset\\
 &\subset &\7q^{(k-1)}+\sigma\7q\;,
\end{eqnarray*}
and the proof is complete.
\qed
\setlength{\abovedisplayskip}{9pt minus 2pt}%
\setlength{\belowdisplayskip}{9pt minus 2pt}%
We are now in the position to characterize holomorphic (non)degeneracy in terms of a
 purely algebraic condition on CR-algebras. 
As already mentioned, a homogeneous CR-germ $(M,o)$ 
is holomorphically nondegenerate if and only it is 
$k$\2-nondegenerate for some 
finite $k.$ 
This follows from  Theorem 11.5.1 in \cite{BER}, applied to the homogeneous  case.
\mn
\begin{theorem}
\label{LieOdege}
\setlength{\parskip}{-6pt}
Let $(\7g,\7q)$ be a given
CR-algebra and $(M,o)$ the corresponding homogeneous CR-germ, 
generically embedded into the germ $(Z,o).$
Let $\7q^{(\bullet)}$ be the filtration by subalgebras as in \ref{ffilt}.ii.
Then , for every integer $k\ge 1$
\begin{enumerate}
\setlength{\itemsep}{-0pt}
\setlength{\itemindent}{-0.1cm}
\item[(i)] $(M,o)$ is $k$\2--nondegenerate if and only if 
$\7q^{(k-1)}\ne\7q^{(k)}=\7q^{(\infty)} ,$ \ and then
$k\le \dim\7q^{(0)}-\dim\7q^{(\infty)}.$ 

\item[(ii)] $(M,o)$ is holomorphically degenerate 
if and only if there exists a complex subalgebra 
$\7r\subset \7l=\gc$ with 
$\7q\subsetneq \7r \subset \7q{+}\sigma\7q.$ The latter condition
implies the existence of a locally equivariant CR-morphism $\Psi:M\to M',$ 
whose
fibres 
are positive-dimensional complex submanifolds of $Z.$ Here, $(M',o')$ is the CR-germ associated with $(\7g,\7r).$
\end{enumerate}
\end{theorem}
\Proof
The first part is an immediate consequence of Proposition \ref{knondeg} and  Observation \ref{Fkato}.
For the proof of the second part of the theorem 
recall that the holomorphic degeneracy of $(M,o)$ is equivalent to the fact that $(M,o)$ 
is not finitely nondegenerate. 
Thus, according to (i)  and Lemma \ref{ffilt}.ii,
there exists $n\in \N$ such that $\7q^{(n)}=\7q^{(n+1)}\ne \7q^{(\infty)}.$
This implies  $[\7q^{(n)},\sigma\7q]\subset \7q^{(n)}+\sigma\7q.$
Since  $\7q^{(n)}$ is a subalgebra by Lemma \ref{ffilt}, $\7q^{(n)}+\sigma\7q$ is a
subalgebra, as well. Define

\mn
\centerline{
$\7r:=\sigma\7q^{(n)}+\7q=\sigma(\7q^{(n)}+\sigma\7q)$
}

\mn
and note that $\7r\supsetneq \7r$ and  $\7r+\sigma\7r=\7q+\sigma\7q.$ 
This proves the existence of $\7r$ as claimed. \hb
Let $(M',o')\subset (Z',o')$
be the CR-germ, associated with the CR algebra $(\7g,\7r).$ 
Since $\7q\subset \7r,$ the identity map on $\7l=\gc$ induces a morphism 
$(\7g,\7q)\to (\7g,\7r),$ and by  Proposition \ref{IsoOfCat}
a CR-morphism $\Psi:(M,o)\to (M',o')$ 
which is the restriction of a holomorphic surjective morphism $\widehat\Psi:(Z,o)\to (Z',o').$ We claim that the germ of the fibre $\Psi^{-1}(o')$
coincides with the germ of the fibre $\widehat\Psi{}^{-1}(o'):$
This can be seen by comparing the dimensions: a simple check shows
that the injection 

\sn
\centerline{$T_o(\Psi^{-1}(o'))=\qu{\7g\cap\7r}{\7g\cap\7q} \longrightarrow
\qu{\7r}{\7q}=T_o(\widehat\Psi{}^{-1}(o'))$}

\sn is also surjective.  \qed

\noindent
{\bf Remarks.}
\break
$\bullet$ In \cite{MN}, certain purely algebraic nondegeneracy conditions of 
CR-algebras have been introduced. 
However, their geometric interpretation, in particular the characterization of
holomorphic (non)degeneracy as  given in the remark 
following Prop.~13.3, compare also Theorem 3.2 in \cite{AM},
 contradicts our Theorem \ref{LieOdege}.

\sn
$\bullet$ A \nn homogeneous, holomorphically degenerate CR-manifold $M$ is,
at generic points in sense of \cite{BRZ}, locally CR-equivalent to a product
of a lower-dimensional CR-manifold and a complex manifold.
This is a consequence of Proposition 3.1. in \cite{BRZ}. 

\bbn
{\bf Minimality in terms of CR-algebras.}
Recall that a  CR-manifold $(M,\H,J)$ is called 
{\it minimal}\/ at $o\in M$
if for each locally closed 
submanifold $Y\subset M$ such that $o\in Y$ and $\H_y\subset T_yY$ for $y\in Y$ 
the identity  $T_oY=T_oM$ holds, i.e., $(M,o)=(Y,o).$ 
In the locally homogeneous situation the property of being minimal at one particular point
is equivalent to the minimality at all points of   $M.$
As before, $\7H=(\7q+\sigma\7q)^\sigma\subset \7g$ and $\H_o\cong\qu{\7H}{\go}.$

\begin{theorem}
\label{LieMini}
Given a CR-algebra $(\7g,\7q),$ let $(M,o)$ be the underlying CR-germ.  
Then
$M$ is minimal  at $o$ if and only if the smallest subalgebra of $\7g,$
which contains $\7H,$ is $\7g$ itself.
\end{theorem} 

\Proof The minimality condition can be reformulated as follows:
Define inductively the following ascending chain of
subbundles (associated with the locally homogeneous CR-manifold $M$):
$$ \H^{(0)}:=\H,\quad \H^{(j)}:= \H^{(j-1)}+[\H^{(j-1)},\H^{(j-1)}] \text\quad \text{for }j>0\;. $$
Here, $[\H^{(\ell)},\H^{(\ell)}]$ stands for the subbundle generated by
all brackets $[\xi,\eta],$ where $\xi,\eta$ run through 
 local sections in $\H^{(\ell)}.$
The minimality of $M$ is equivalent to the
condition $\bigcup_{k\ge 0} \H^{(k)}=TM.$ 
In our situation all subbundles $\H^{(k)}$ are homogeneous; 
hence, they are completely determined by the fibres at $o\in M.$
Let $\7H^{(k)}\subset \7g$ denote the subspaces containing  $\go$
such that $\H^{(k)}_o=\qu{\7H^{(k)}}{\go}$ for all $k.$
Observe that the map

\sn
\hfil $\Gamma(M,\H^{(k)}) \times \Gamma(M,\H^{(k)})\longrightarrow 
\qu{\Gamma(M,\H^{(k+1)})}{\Gamma(M,\H^{(k+1)})},$ 

\sn
given by the Lie brackets is $C^\infty(M)$\2--bilinear. Consequently,
 we can employ the Main Lemma \ref{bas}: The corresponding tensor $\H^{(k)}_o\times \H^{(k)}_o \to \qu{\H^{(k+1)}_o\!}{\,\H^{(k)}_o}$ 
is simply given by the Lie bracket in $\7g.$ This yields an inductive definition of all $\7H^{(k)}:$
The subspace $\7H^{(k+1)}$ is generated by elements $u\in \7H^{(k)}$
and all Lie brackets $[u,v]_\7g,$ $u,v\in \7H^{(k)}.$ 
If the smallest Lie algebra in $\7g$ which contains $\7H^{(0)},$ 
coincides with  $\7g$ then $\bigcup_{k\ge 0} \7H^{(k)}=\7g$ and consequently
$\bigcup_{k\ge 0} \H^{(k)}=TM,$ i.e., $M$ is minimal. The opposite direction,
i.e., ``$M$ minimal implies $\7g$ is the smallest subalgebra containing $\7H$''
is easier to see: The existence of a  proper subalgebra
of $\7g$ which contains $\7H,$ would imply the existence of an integral 
manifold (Nagano leaf) through $o,$ strictly lower-dimensional than $M.$ 
But this contradicts the minimality of $M.$\qed

\bn
{\bf Orbits in flag manifolds.}
In this subsection let  $Z$ stand for a flag manifold, i.e., a projective homogeneous manifold with $b_1(Z)=0.$  Let
$L\subset \mathrm{Aut}_\5O(Z)$ be a complex subgroup which
acts transitively on $Z,$ i.e., $Z=L/Q.$
In such a case $L$ is semisimple and the isotropy subgroup $Q$ is
 parabolic. Select a real form
$G$ of $L.$ The $G$\2-orbits in $Z$ provide a broad class of
examples of CR-manifolds. 
For instance, all bounded symmetric domains $D\subset \C^N$ and the pieces of the natural stratification of their boundaries
arise as certain  orbits of the above type.
In \cite{KZ2}, global actions  of so-called real forms of tube type 
have been considered in the particular case  where
 $Z$ is a Hermitian compact symmetric space.
Recall that a real form $G$ of a complex semisimple Lie group $L$ is called
of {\it tube type}\/ if $G$ has an open orbit in $Z,$ 
which is biholomorphically equivalent to a bounded symmetric domain of 
tube type.
It has been proven 
(Theorem 4.7 in \cite{KZ2}) with Jordan algebraic tools
that for an $G\subset L$ of tube type   
 each $G$\2--orbit $M$ in  an irreducible Hermitian space 
$Z=L/Q,$ 
which in neither open nor totally real 
is 2-nondegenerate and minimal.
As shown in \cite{W} in each flag manifold $Z$ there is precisely one closed $G$\2--orbit $Y.$
Further, $\dim_\R Y\ge \dim _\C Z,$ and the closed orbit is totally real if and only if
  $\dim_\R Y= \dim _\C Z,$ as  is the case for $G$ of tube type and $Z$ the corresponding
Hermitian space.

\sn

Natural generalizations of irreducible compact
Hermitian symmetric spaces are  the flag manifolds
$Z=L/Q$ with  second Betti number $b_2(Z)$ equal to 1,  
or equivalently, where $Q$ is  a 
maximal parabolic subgroup.
In this situation  Theorem 4.7 from \cite{KZ2} can be generalized as follows:

\begin{theorem}
\label{kaza} Let $L$ be a complex simple Lie group, $G\subset L$ an arbitrary
 real form and $Q\subset L$ a parabolic subgroup. 
\setlength{\parskip}{-4pt}
\begin{enumerate}
\setlength{\itemsep}{-1pt}
\setlength{\itemindent}{0.0cm}
\item[\rm(\it i\rm)]
Assume that $Q$ is\/ {\bit maximal} parabolic.
Then every non-open $G$\2--orbit $M$ in $Z:=L/Q$ is 
holomorphically nondegenerate. All such orbits are also  minimal, 
except for the totally real ones.
\item[\rm(\it ii\rm)] In particular, if
 $Z=L/Q$ is an irreducible Hermitian space with $L=\mathrm{Aut}(Z)^\circ$ and $G\subset L$ an 
{\em arbitrary} real form then every $G$\2--orbit $M$ which is not open 
is $k$\2--nondegenerate with $k\le 2.$  For every
such orbit $M,$ which in addition is not totally real,
 $(M,Z)$ belong to the class $\7C$ in the sense of {\rm 4.4.}~in \rm\cite{KZ2}.\it 
\item[\rm(\it iii\rm)] If $Q$ is not maximal, then there always
exist  non-open holomorphically\/ {\bit degenerate} $G$\2--orbits in $Z.$

\end{enumerate}
\end{theorem}
\Proof Let $\sigma:\7l\to \7l$ be the involution given by the real form 
$G\subset L.$ Let $\7q_z$ be the isotropy Lie algebra at a point $z\in Z=G/Q,$
$M:=G\cd z$ the orbit with the inherited CR-structure such that 
neither $\7q_z+\sigma \7q_z=\7g$ (i.e., $M$ is not open) nor $\7q_z+\sigma\7q_z=\7q$ (i.e., $M$ is not totally real; here we follow the notational
convention from \cite{W}
and denote the complex isotropy at $z\in Z$ by $\7q_z$ rather than $\7l_z$).
Since the only Lie algebra, properly containing $\7q_z$  (and in particular 
$\7q_z+\sigma\7q_z$),  is $\7l$ itself, Theorem
\ref{LieOdege} together with
 Theorem \ref{LieMini} imply the first part of the claim.
The estimate for the order of nondegeneracy $k$ in the Hermitian case $Z=L/Q$
follows from Theorem \ref{hypersurfOrb} together with the following well-known
technical fact that $c(\7q)=1$ (\cite{W1}, 
see our notation in the paragraph preceeding Theorem \ref{hypersurfOrb}).
As the example $\P_{\!2n-1}=\SL_{2n}(\C)/Q= \Sp_n(\C)/P$ shows, 
complex Lie groups of different dimensions 
may act transitively on a given flag manifold.
\hb
If $Q$ is not maximal, there exists a maximal parabolic $Q',$ 
containing $Q,$ such that $G$ is not transitive on $L/Q'.$ Further, 
there is the $L$\2--equivariant holomorphic map 
$ \pi:\qu LQ\to\qu L{Q'}=:Z'$ with positive-dimensional complex fibres.
Let $M'\subset Z'$ be an arbitrary $G$\2--orbit which is not open.
Then $\pi^{-1}(M')$ consists of finitely many $G$\2--orbits.
In particular there exists an orbit $M$ which is open in $\pi^{-1}(M').$
The fibres of the restriction $\pi:M\to M'$ (which is a CR-map) are then 
complex manifolds and consequently $M$ is locally equivalent to
a product of a CR-manifold and a positive dimensional complex manifold. 
This implies that $M$ is holomorphically degenerate.
\qed

\section{\llap{.\kern.75em}A 3-nondegenerate homogeneous CR-manifold}
\setcounter{equation}{0}

The purpose of this section is to give an explicit example of a homogeneous 
$3$\2-nondegenerate CR-manifold. Recall that all 
CR-manifolds which occur  in \cite{KZ2} are either holomorphically degenerate
or 2-nondegenerate. In the Hermitian symmetric spaces, complementary
to those considered in  \cite{KZ2} there are also 1-nondegenerate CR-manifolds.
Up to our knowledge, there are no known examples of {\it homogeneous}\/
$k$\2-nondegenerate CR-manifolds with $k\ge 3.$ 

\sn

Examples, promising
in search of homogeneous CR-manifolds with higher nondegeneracy, arise as
 orbits of real forms in flag manifolds.
Note however that the Jordan-algebraic methods used
in \cite{KZ2} in the particular case, where $Z$ is
 a compact Hermitian symmetric space, cannot be generalized to
arbitrary flag manifolds. Nevertheless, this (bigger) class of 
orbits of real forms with induced CR-structures coming from
general flags $L/Q$
is still quite accessible from a computational point of view:
This is due to the fact that
every complex isotropy Lie algebra $\7q=\7l_z=:\7q_z$ contains a 
$\sigma$\2--stable Cartan subalgebra $\7t$ (\cite{W}, Thm.~2.6).
Here, $\sigma:\7l \to \7l$
is the conjugation induced by the real form $\7g\subset \7l.$  
Consequently, all 
subspaces $\7q^{(\ell)}$ of $\7l$
 contain this Cartan subalgebra
and are direct sums of root spaces. The algebraic manipulation
of the corresponding CR-algebra boil down to the combinatorics 
of root systems. In the next subsection we explain for a  particular example
all that in greater detail.

\bn
{\bf The example, described in geometric terms.}
In the context of flag manifolds, 
the simplest example of a 3-nondegenerate CR-manifold
arises as a (locally closed) hypersurface orbit $\4M:=G\cd z_0$ in $Z:=\G_2^b(V).$
Here,  $V$  is isomorphic to $\C^7;$ further, 
\hb
$\bullet$ $b:V\times V\to \C$  is  a symmetric nondegenerate 2-form;
 it determines the orthogonal group $L:=\mathrm{Isom}(V,b)\cap \SL(V)=\Iso(V,b)^\circ\cong \SO(7,\C);$ and we write  $\Iso(V,b)=\{A\in \GL(V):b(Av,Aw)=b(v,w)$ for all $v,w\in V\},$
\hb
$\bullet$ $\G_2^b(V)$  denotes the Grassmannian of $b$\2--isotropic
 2--planes in $V;$ it is a 7-dimensional submanifold
of the Grassmannian $\G_2(V).$ Further, $\mathrm{Aut}_\5O(Z)^\circ=L.$
\hb 
$\bullet$ $G\subset \mathrm{Aut}_\5O(Z)^\circ=L$ is a real form,
which is isomorphic to $\SO(3,4)^\circ.$ Such a $G$
is determined by
an appropriate choice of a  maximal and totally real subspace 
$\R^7\cong V^\R\subset V:$ 
we have  $G=\{g\in L: g(V^\R)=V^\R\}^\circ\cong\SO(3,4)^\circ.$
(A similar construction remains valid for all real forms of type $\SO(p,q)$ in
 $\SO(p{+}q,\C).$)  
Write $\underline \sigma:V\to V$ for the anti-linear conjugation with $V^{\underline\sigma}=V^\R$ and $\sigma:L\to L\subset \GL(V)$ for $\sigma(g)=\underline \sigma\circ g \circ \underline \sigma.$ For example, if  $V:=\C^7,$ \  $b(z,w):=\sum_{j=1}^7 z_jw_{8-j},$ \ and 
$L:=\{A\in \C^{7\times 7}:b(Ax,Ay)=b(x,y)\ \ \forall\ x,y\in \C^7\}^\circ,$ then we have 

\mn
\Quad\qqquad $G=L\cap \R^{7\times 7}$ \quad and \quad  $\sigma(A)=\overline A.$

\mn
$\bullet$  Finally, define the associated Hermitian 2-form $h^b(v,w):=b(v,\underline\sigma(w)).$ It has  signature $(3,4)$
and  $G=(\Iso  (V,b)\cap \Iso(V,h^b))^\circ.$   

\bn

Let $\4H\subset \G_2^b(V)$ be the set of all planes $E\in \G_2^b(V)$ such that 
$h^b\rest{E}$ is degenerate. This  is  a (singular) real hypersurface in $Z,$ stable
under $G.$ 
The CR-manifold $\4M$ is  a $G$\2--orbit in the smooth part of $\4H,$ open in $\4H.$ 
Note  that 
the closed $G$\2--orbit $\4Y$ in $\G_2^b(V)$ is totally real and isomorphic to the
real Grassmannian $\G_2^b(V^\R).$ It is also contained in $\4H.$

\mn 

The geometrically described hypersurfaces $\4H$ and $\4M$
can also be given in local coordinates as the zero set of a  function $\rho.$ 
Note that $Z$ as  a (7-dimensional) flag manifold is covered by Zariski open subsets 
$U\subset Z$ which are all isomorphic to $\C^7$ and provide coordinate charts on $Z.$ 
We pick one of such charts,
$U,$ centered in a point of the totally real orbit $\4Y\subset Z,$  and  give a
 defining function for $U\cap \4M.$ 
Write  $\text{\9z}:=(z_1,z_2,z_3)$ and $\text{\9w}:=(w_1,w_2,w_3)$
for (row) vectors in $\C^3,$ define the quadratic 2-form
$c(\text{\9z},\text{\9w}):=\frac12  (z_1w_3+z_2w_2+z_3w_1)\;, $
and write $(\text{\9z,\ \9w},u)$ for (row) vectors in $\C^7\cong U.$
Then the function $\rho=\rho(\text{\9z},\text{\9w},u)$ is polynomial of degree 4
and is
given as the determinant of a $2\times 2$ matrix:
\begin{equation} 
\rho=\det{\small
\begin{pmatrix}c(\text{\3w},\text{\3w})-2c(\overline{\text{\3w}}, \text{\3w})+c(\overline{\text{\3w}},\overline{\text{\3w}}) &
\kern-4pt c(\text{\3z},\text{\3w})-2c(\overline{\text{\3z}}, \text{\3w})+c(\overline{\text{\3z}},\overline{\text{\3w}}) +i\Im u \\ \\
c(\text{\3z},\text{\3w})-2c({\text{\3z}}, \overline{\text{\3w}})+c(\overline{\text{\3z}},\overline{\text{\3w}}) -i\Im u &c(\text{\3z},\text{\3z})-2c(\overline{\text{\3z}}, \text{\3z})+c(\overline{\text{\3z}},\overline{\text{\3z}})
\end{pmatrix}\;.
}
\end{equation} 
Note that \wrt our coordinates chart, $\R^7=U\cap \4Y,$ i.e., $\R^7\cap \4M=\nix.$ 
However, $(0,i,0,0,0,0,0)\in \4M.$

\mn

Instead of a direct examination of this equation (which might be
 one possibility
to check that $\4M$ is uniformly 3-nondegenerate), we give a description
of the corresponding CR-algebra and use Theorem \ref{LieOdege} 
to check that order of nondegeneracy.
The method used below can actually be generalized to find the 
associated CR-algebra of
an arbitrary
 $G$\2--orbit in an arbitrary flag manifold $L/Q.$ For simplicity, 
we restrict our considerations to the particular case of our 
hypersurface orbit in $\G_2^b(V)=L/Q.$

\bn

\mn
{\bf A root theoretical description and further generalization.}
Our first task is to identify the conjugacy class of the parabolic isotropy 
subalgebra $\7q\subset \7l$ (of $\G_2^b(V)$) in terms of root subsystems.
For the general theory of parabolics we refer to \cite{Ser}.
 Recall  that every parabolic subalgebra $\7q$
contains a Borel subalgebra $\7b$ (a maximal solvable subalgebra of $\7l$)
and a Cartan subalgebra $\7t$ 
(a maximal subalgebra, containing semisimple elements only) such that
$\7t\subset \7b\subset \7q.$ 
The Lie algebra $\7l$ has a decomposition
$ \7l=\7t\oplus \bigoplus_{\alpha\in\Phi(\7l,\7t)}\7l^\alpha,$ $\7t=\7l^0,$ 
into the root spaces $\7l^\alpha=\{v\in \7l:[d,v]=\alpha(d)\cd v$ for all $d\in \7t\}$ 
\wrt the Cartan subalgebra $\7t.$
Here, $\Phi:= \Phi(\7l,\7t)\subset \7t^*$ stands for 
the set of roots (i.e., non-trivial eigen-functionals $\alpha:\7t\to \C,$
 which appear in the root decomposition).
It is possible to select a  $\sigma$\2--stable  Cartan subalgebra $\7t\subset \ ( \7b \ \subset \ ) \  \7q.$ In such a case $\sigma$ induces a permutation $\sigma:\Phi\to \Phi$ of roots.
We follow here
the general
convention and declare $\Phi(\7b,\7t)$ to be the negative roots $\Phi^-.$ 
Let $\Pi\subset \Phi^+$ denote the corresponding simple roots. The  conjugacy classes of parabolic subalgebras of $\7l$ are parameterized 1-to-1 by
 subsets of $\Pi:$ This assignment is given by
 $\7q\leftrightsquigarrow \4Q^r:=\Phi(\7q,\7t)\cap \Pi.$ 
(Some authors use the complementary identification
$\7q\rightsquigarrow \Pi\smallsetminus (\Phi(\7q,\7t)\cap \Pi).$) 

\mn

In our particular example we have $L\cong \SO_7(\C),$ and 
{\large $\hbox{\mathsurround=0pt$\circ$}
\vrule width17pt height3.2pt depth-2.7pt 
\hbox{\mathsurround=0pt$\circ$} 
\vrule width19pt height4.3pt depth-3.8pt\kern-19pt
\vrule width19pt height2.3pt depth-1.7pt
\kern-17pt>\hbox{\mathsurround=1pt$\circ$}$} is the Dynkin diagram of $\Pi.$ 
Let
$\alpha_1,\alpha_2,\alpha_3$ denote the consecutive simple roots, with $\alpha_3$ short. 
Then the parabolic isotropy subalgebra, defining
$\G_2^b(V),$ corresponds to the subset $\4Q^r:=\{\alpha_1,\alpha_3\}.$
Let $\7q=\7q_z$ be the complex isotropy at $z$ in a  given $G$\2--orbit.
In our case, the computation of the various
subspaces $\7q^{(\infty)}\subset \7q^{(\ell)}\subset \7q$ (see \ref{ffk}) and
$\7q+\sigma\7q$ can be reduced to the computation of the corresponding subsets of $\Phi$
which, in turn, is pure combinatorial:
Select  a $\sigma$\2--stable Cartan subalgebra and a Borel subalgebra with
$\7t\subset \7b\subset \7q_z.$
In our particular example $\4M,$ the induced action of $\sigma$ on the roots $\Phi=\Phi(\7l,\7t)$ is depicted in the
figure below:  \kern3pt For short,\break

\parskip=-8pt 
\noindent
{\parbox[b]{9.4cm}{the digits 
stand for
 the coefficients in the expression of a root $\beta$ \wrt the basis $\alpha_1,\alpha_1,\alpha_3.$ For instance, ``$-012$'' stands for $-\alpha_2-2\alpha_3$ and 
$\7l^{-012}:=\7l^{-\alpha_2-2\alpha_3}.$ The arcs connect all pairs $\beta,\sigma(\beta),$
hence, completely determine $\sigma:\Phi\to \Phi.$ A glance at that diagram immediately shows  that
\abovedisplayskip=9pt
\belowdisplayskip=0pt
$$ \diagram{
\7q^{(\infty)}=\7q^{(3)} & \ = \ & \7t\oplus \7l^{100}\oplus \7l^{-112}\oplus \7l^{-001}\oplus \7l^{-011}\oplus \7l^{-012}
\cr \hfill \7q^{(2)}&=& \7q^{(3)} \oplus  \7l^{-122} \hfill\cr 
\hfill\7q^{(1)}&=& \7q^{(2)}\oplus \7l^{-010}\oplus \7l^{-111}\hfill
\cr
\hfill\7q=\7q^{(0)}&=& \7q^{(1)}\oplus \7l^{-100}\oplus\7l^{-110}\oplus \7l^{001}\;\;.\hfill }
$$
}}\hfill\epsfig{file=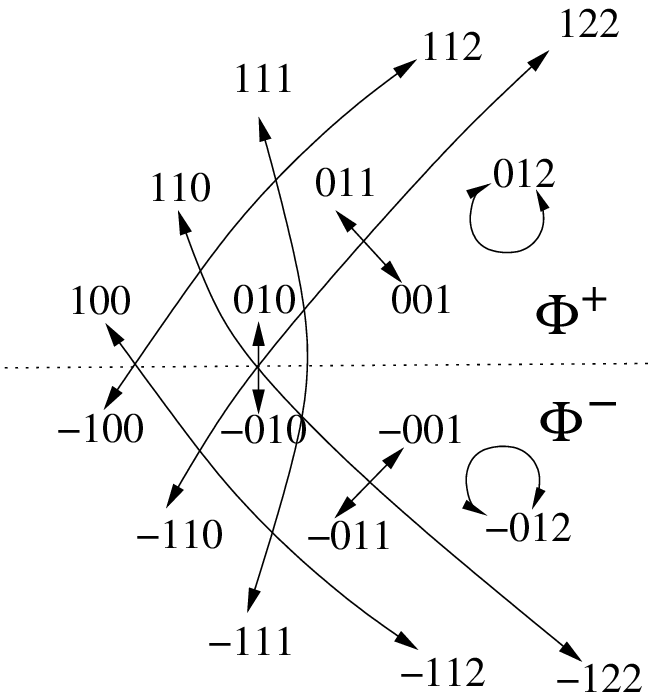,height=4.7cm}
 
\parskip=5pt\noindent
The particular shape of $\sigma:\Phi\to \Phi$ could be computed by ``brute force''
simply by selecting a base point $z\in \4M,$ describing the corresponding
subalgebras $\7t=\sigma(\7t)\subset \7b\subset \7q=\7q_z$ in terms of $7\times 7$ matrices and 
finally computing the induced involution $\sigma:\Phi(\7t)\to \Phi(\7t).$ 
A more elegant way, suitable for a generalization to arbitrary orbits in flag manifolds is the following: Start with a point $y\in \4Y$ on the closed orbit.
In our example, for $\7t'\subset \7b'\subset \7q'=\7q_y,$ the action of $\sigma$ on $\Phi'=\Phi(\7t')$ 
is particularly simple: It is the identity (in the general case $\sigma:\Phi'\to \Phi'$ 
can be read off the Satake diagram of the real form $\7g$).  
Apply certain partial Cayley transformations
   to obtain
$\7t:=\cc(\7t')\subset \cc(\7b')\subset \cc(\7q_y)=\7q_z$ and $\cc(y)=z$ such that $z$ is contained in the orbit  under consideration.
In our particular case, we perform the partial Cayley transformations \wrt
the strongly orthogonal roots $\gamma_1:=\alpha_1+\alpha_2+\alpha_3$ and 
$\gamma_2:=\alpha_2$ in $\Phi',$ i.e., $\cc=\cc_{\gamma_1}\circ \cc_{\gamma_2}.$
The induced involution $\sigma$ on $\Phi(\7t)=\cc(\Phi')$ can be computed 
by a repeated application of  the  formula
\[
\sigma \cc_\gamma(\beta)=\cc_\gamma(\sigma\beta)-\langle \beta\,|\,\gamma\rangle\cd \cc_\gamma(\gamma)  \qquad \beta,\gamma\in \Phi'
\;,
\]
where $\kappa$ is the symmetric product on $\Phi',$ induced by the Killing form and $\langle \beta\,|\,\gamma\rangle:=2\,\frac{\kappa(\beta ,\gamma)}{\kappa(\gamma,\gamma)}.$
This method can be used to handle arbitrary   orbits of real forms in arbitrary flag manifolds.

\parskip=0pt
\mn

Let $x\in \4M$ be a point and $\hol(\4M,x)$ the Lie algebra of germs of 
all infinitesimal CR-transformations at $x$ (see Section 3).
By the nondegeneracy of $\4M$  we have $\dim \hol(\4M,x)<\infty,$
and clearly  $\so(3,4)=\7g\subset \hol(\4M,x).$ 
We do not know, however,  if this inclusion is proper.
It should be noted that Prop.~3.8 in \cite{KZ2} which uses the existence of  nonresonant vector fields
 does not apply
in our situation since  $\G_2^b(\C^7)$
is not a Hermitian space: Due to the following lemma  there is no nonresonant vector field on $\G_2^b(\C^7)$
coming from $\7l=\so(7,\C).$ 
More precisely,

\begin{lemma} 
\label{} Let $M=G\cdot z\subset Z=L/Q$ be an orbit of a real form in an arbitrary 
 flag manifold such that $\7l$ is simple. Then $\7l$ contains a nonresonant
vector field if and only if $Z$ is Hermitian and $\7l=\aut(Z).$ In the 
Hermitian case there always exists a nonresonant vector field in 
$\aut(Z)\subset \hol(Z,z)\cong\hol(\C^N,0)$
with linear part equal to $\Id.$
\end{lemma}
\Proof It is sufficient to consider the isotropy action of a Cartan subalgebra 
$\7t\subset \7q$ on $T_zZ\cong \7q^n.$ This action  is diagonalizable and the
corresponding eigenfunctionals $\beta\in\7t^*$ (i.e., roots in $\Phi^n$)
 determine the eigenvalues of the linear parts of the vector fields, 
induced by elements in $\7t.$  We use here the  decomposition
 $\Phi=\Phi^{-n}\cup \Phi^r\cup \Phi^n$ of the root system $\Phi=\Phi(\7l,\7t),$ 
induced by $\7q$ such that 
 $\Phi(\7q,\7t)=\Phi^{-n}\cup \Phi^r$ are the roots of the nilpotent resp.
reductive part of $\7q.$ 
If $\7q$ is {\it not}\/ of Hermitian type (i.e, $Z=L/Q$ is not a Hermitian symmetric 
space with $\7l=\aut(Z)$)
then there always exist $\alpha,\beta\in \Phi^n$ with $\alpha+\beta\in \Phi^n.$ This 
violates the nonresonance condition (as given in \cite{KZ2}).  
The Hermitian situation is well-known.\qed
\noindent

\mn

The above defined hypersurface $G$\2-orbit $\4M$ is a particular example of a finitely nondegenerate CR-manifold. 
One would expect that there are
$G$\2-orbits in flag manifolds
which are finitely nondegenerate of  arbitrary high order. Surprisingly,
at least for hypersurface orbits,
 this is not true as the following theorem shows.
Before we state it, we recall some standard notation:
Given a parabolic subalgebra $\7q\subset \7l,$ select $\7t\subset \7b\subset
\7q,$ ($\7t$ a Cartan and $\7b=\7t\oplus \bigoplus_{\Phi^-}\7l^\gamma $ a Borel subalgebra),
and let $\Pi=\{\alpha_1,...\,,\alpha_q,...\,,\alpha_n\}
\subset -\Phi(\7b,\7t)$ be the corresponding simple roots.
If  $\7q$ is maximal, it is determined by a subset $\Pi\sm \{\alpha_q\}=\Phi(\7q,\7t)\cap \Pi,$ where 
$\alpha_q$ is a simple root. 
Let $c(\7q):=\max\{ c_q(\gamma):\gamma\in \Phi^+\},$
where $c_q(\gamma)=c_q$ is the $q^\text{th}$ coefficient in the expression
$\gamma=c_1\alpha_1+\cdots +c_q\alpha_q+\cdots +c_n\alpha_n.$
For example $c(\7q)=1$ if $Z=L/Q$ is a Hermitian space with $\mathrm{Aut}(Z)^\circ=L$ (see \cite{W1}).

If $L=\prod L_j$ is a direct product of simple complex Lie groups and
$G=\prod G_j$ a real form of $L$ such that $G_j$'s are arbitrary
real forms in the simple factors $L_j$'s, the corresponding $G$\2--orbit $M$
in $Z=L/Q=\prod L_j/Q_j$ is also a direct product $M=\prod M_j$ 
as a CR-manifold.
Consequently, we may restrict our consideration to the case $L/Q$ where $L$
is  simple:

\begin{theorem}
\label{hypersurfOrb}
Let $Z:=L/Q$ be an arbitrary flag manifold where $L$ is a simple complex group
and  $G\subset L$ a real form.
Let $M:=G\cd z$ be an orbit in $Z.$
\setlength{\parskip}{-4pt}
\begin{enumerate}
\setlength{\itemsep}{-1pt}
\setlength{\itemindent}{0.0cm}
\item[$(i)$] Assume that $M$ is a real hypersurface in $Z.$
Then  $\!M$ is holomorphically nondegenerate
 if and only if  $Q$ is a {\em maximal} parabolic subgroup of $L.$
\item[$(ii)$]
Assume that $b_2(Z)=1,$ i.e., $\7q$ is maximal parabolic, and $M$ is not open in $Z.$
Let $k(M)$ denote the order of nondegeneracy of $M.$
Then $k(M)\le c(\7q)+1\le 7$ (with  $c(\7q)$ as defined above).
In particular, $k(M)\le 3$ if $L$ is not an
exceptional simple group.
\end{enumerate}
\end{theorem} 
\Proof
Let $M=G\cd  z\subset L/Q$ be a hypersurface orbit and $Q=Q_z$ the complex 
isotropy subgroup at $z.$ As explained before,
select   $\7t\subset \7b\subset \7q$ with a $\sigma$\2--stable Cartan subalgebra
$\7t.$
The assumption $\text{codim}_Z(M)=1$ implies that $\7q+\sigma\7q$ is a hyperplane in $\7l,$ i.e., there exists $\gamma\in \Phi(\7t)$ with $\sigma(\gamma)=\gamma$
such that $\7l=\7q{+}\sigma\7q \ \oplus \7l^\gamma.$ 
Let $\4Q^r\subset \Pi=\{\alpha,\ldots,\alpha_n\}$ 
be the subset determined by $\7q$ and $\gamma=\sum n_j\alpha_j.$ Write $\mathrm{supp}_\Pi(\gamma):=\{\alpha_j\in \Pi: n_j>0\}.$
Clearly
$\mathrm{supp}_\Pi(\gamma)\not\subset \4Q^r.$ Select an element $\hat\beta\in \mathrm{supp}_\Pi(\gamma)\smallsetminus \4Q^r.$ 
Then $\Pi\sm\{\hat\beta\} \ \supset \4Q^r$ and consequently the parabolic subalgebra $\widehat {\7q}$ corresponding to the set $\Pi\sm\{\hat\beta\}$ contains $\7q$
and we have $\7q\subset \widehat{\7q}\subset \7q+\sigma\7q.$
If $\widehat{\7q}\supsetneq \7q,$ i.e., $\7q$ is not a maximal parabolic,
the corresponding orbit $G\cdot z$ is holomorphically degenerate, due to Theorem \ref{LieOdege}. This proves the first part of the statement. 
\hb
To prove the second part, 
let $\7q=\7q_x$ be the complex isotropy subalgebra at $x\in M\subset Z$ and
$\7l=(\7q+\sigma\7q)\oplus \7l^\Gamma$ where $\7l^\Gamma:=\bigoplus_{\gamma\in \Gamma} \7l^\gamma,$
i.e., $|\Gamma|$ is the CR-codimension of $M$ in $Z.$
Let $\alpha_q\in \Pi$ be the simple root such that
$\Phi(\7q,\7t)\cap \Pi=\Pi\sm\{\alpha_q\}.$
Note that then $\alpha_q\in \mathrm{supp}_\Pi(\gamma)$ for every $\gamma\in \Gamma.$ 
Further, the simple root $\alpha_q$ determines the following $\Z$\2--filtration 
$\bigoplus_{-\infty}^\infty \7l_j,$ where the homogeneous parts are given by
\[ 
\7l_j:=
\begin{cases}
\bigoplus_{c_Q(\beta)=j} \7l^\beta & \text{for } j\ne0 \quad \  (\beta=c_1(\beta)\alpha_1+\cdots +c_q(\beta)\alpha_q+\cdots +c_n(\beta)\alpha_n)  
\\
\7q^\text{red} & \text{for } j=0 \quad \ \text{\small(Here, $\7q^\text{red}$ is the reductive part of $\7q,$ containing $\7t$)}
\end{cases}
\]
and $\7l_j=0$ for $|j|>c:= c(\7q).$
We have then $\7q=\bigoplus_{j=-c}^0\7l_j.$ For short, write 
\begin{center}
$\displaystyle\7q^{n\sm\Gamma}=\bigoplus_{\textstyle{c_q(\beta)>0\atop \beta\not\in\Gamma}} \kern-.5em\7l^\beta \quad \ \subset\ \7q^n$
\end{center}
and note that $\sigma\7q=\7q^{n\sm\Gamma}\oplus (\7q\cap\sigma\7q) =
\7q^{n\sm\Gamma}\oplus\7q^{(\infty)}.$ Let $\7q=\7q^{(0)}\supset \7q^{(1)}\supset \cdots \supset \7q^{(\infty)}$ be the filtration as defined in \ref{ffk},
and write $\7q^{(\ell)}_{-c}\oplus \cdots \oplus \7q^{(\ell)}_{-1}\oplus
\7q^{(\ell)}_{0}$ for the corresponding gradation of the $\7q^{(\ell)}$'s, $\ell=0,1,2...$ (and $\7q^{(\ell)}_j=\7q^{(\ell)}\cap \7l_j$).
Since $\7t\subset \7q^{(\infty)}$ and  all
 subalgebras $\7q^{(\ell)}$ are $\ad(\7q^{(\infty)})$\2--stable, the 
condition defining $\7q^{(\ell)}$ (see \ref{ffk}) is equivalent to 

\mn
$(\star)$\qqquad $\7q^{(\ell)}=\bigoplus \7q^{(\ell)}_j$\quad with \quad $[\7q^{(\ell)}_j,\7q^{n\sm \Gamma}]\subset \7q^{(\ell-1)}_{j+1}\oplus \cdots \oplus  \7q^{(\ell-1)}_{0}\oplus \7q^{n\sm\Gamma} \;.$

\mn
The last statement of the theorem follows then from the following

\mn
{\bf Claim.} For every $p\ge 0$ we have $\7q^{(p+1)}_{-p}=\7q^{(p+2)}_{-p}=\cdots=\7q^{(\infty)}_{-p}.$

\mn
We prove the claim by induction, using condition $(\star).$
 For $p=0$ we have $\7q^{(1)}_0=\7q^{(2)}_0=\cdots$ since the condition $[\7q^{(j)}_0,\7q^{n\sm\Gamma}]\subset \7q^{n\sm\Gamma}$ does not depend on $j.$ Assume, we have already proved the claim
for $p=0,1,...\,,m.$
Then, for every $\ell\ge m+2,$ we have
\begin{eqnarray*}
[\7q^{(\ell)}_{-m-1},\7q^{n\sm\Gamma}] &\subset & \7q^{(\ell-1)}_{-m}\oplus \7q^{(\ell-1)}_{-m+1}\oplus\cdots \oplus\7q^{(\ell-1)}_{0}\oplus \7q^{n\sm\Gamma}=
\\
& &=\7q^{(m+1)}_{-m}\oplus \7q^{(m+1)}_{-m+1}\oplus\cdots \oplus\7q^{(m+1)}_{0}\oplus \7q^{n\sm\Gamma}
\end{eqnarray*}
and consequently the conditions imposed  on 
$\7q^{(\ell)}_{-m-1}\subset \7q^{(m+2)}_{-m-1}$  for each   $\ell\ge m+2$ 
do not depend on  $\ell.$ 
This proves the claim.

\mn
Due to the above claim, at most after $c(\7q)+1$ steps
the filtration  $\7q^{(\bullet)}$ becomes stationary. The theorem follows now from this observation and  Theorem \ref{LieOdege}. 
The values of $c(\7q)$ are bounded by the highest coefficient $c(\7l)$
of the highest root of $\7l.$
A glance at the table of highest roots for the classical and exceptional 
simple Lie algebras $\7l$ yields $c(\7q)\le c(\7l)\le 2$ in the classical cases and
$c(\7e_6)=c(\7g_2)=3,$ $c(\7e_7)=c(\7f_4)=4$ and $c(\7e_8)=6$ in the exceptional cases.
\qed

\sn

\vbox{\noindent
{\bf Problems.} Let $M$ stand for a $G$\2-orbit in a flag manifold $Z=L/Q.$  Generalizing the above methods:
\setlength{\parskip}{-7pt}
\begin{enumerate}
\setlength{\itemsep}{-3pt}
\item[(A)] 
Carry out the case where $\mathrm{CR\text{-}codim}_ZM\ge 2$ and $b_2(Z)\ge 2.$
\item[(B)] Carry out the ``group case,'' i.e., describe the degeneracy of the $G$\2--orbits $M$ in $Z=L/Q,$ where
the real form $G$ carries a complex structure, i.e.,  $L\cong G\times G.$  
\end{enumerate}
}

\setlength{\parskip}{-7pt}

\setlength{\parskip}{-1pt}
\begin{center}
\small  Gregor Fels, Mathematisches Institut, Auf der Morgenstelle 10, 72076 T\"ubingen, Germany
\\
\tt gfels@uni-tuebingen.de
\end{center}
\end{document}